\documentclass{amsproc}

\newtheorem{theorem}{Theorem}[section]
\newtheorem{proposition}[theorem]{Proposition}

\newtheorem{corollary}[theorem]{Corollary}

\theoremstyle{definition}
\newtheorem{definition}[theorem]{Definition}
\newtheorem{example}[theorem]{Example}

\theoremstyle{remark}
\newtheorem{remark}[theorem]{Remark}

\numberwithin{equation}{section}

\begin{document}

\title{Ratliff-Rush Monomial Ideals}

\author{Veronica Crispin Qui{\~n}onez}
\address{Department of Mathematics, Stockholm University, SE-106 91
  Stockholm, Sweden}
\email{veronica@math.su.se}

\date{January 1, 1994 and, in revised form, June 22, 1994.}


\begin{abstract}
Let $I$ be a regular $\mathfrak m$-primary ideal in $(R,\mathfrak
m,k)$. Then its Ratliff-Rush associated ideal $\bar I$ is the largest
ideal containing $I$ with the same Hilbert polynomial as $I$. In this paper we
present a method to compute Ratliff-Rush ideals
for a certain class of monomial ideals in the rings $k[x,y]$ and
$k[[x,y]]$.  We find an upper bound for the Ratliff-Rush reduction number
for an ideal in this class. Moreover, we establish some new characterizations of when
all powers of $I$ are Ratliff-Rush.
\end{abstract}

\maketitle

\section{Introduction}

Let $R$ be a Noetherian ring and let an ideal $I$ in it be regular,
that is, let $I$ contain a nonzerodivisor. Then the ideals
$(I^{l+1}:I^l),\; l\geq 1,$ increase with $l$. The union $\tilde
I=\bigcup_{l\geq 1}^{\infty} (I^{l+1}:I^l)$ was first studied by
Ratliff and Rush in \cite{RR}. They show that $(\tilde I)^l=I^l$ for
sufficiently large $l$ and that $\tilde I$ is the largest
ideal with this property. Hence, $\tilde{\tilde I} =\tilde
I$. Moreover, they show that $\widetilde{I^l}=I^l$ for sufficiently
large $l$. We call $\tilde I$ the Ratliff-Rush ideal associated to
$I$, and an ideal such that $\tilde I=I$ a Ratliff-Rush ideal. The
Ratliff-Rush reduction number of $I$ is defined as $r(I)=\min\,\{l\in\mathbb Z_{\geq
  0}\:\vert\:\tilde I=(I^{l+1}:I^l)\}$. 

The operation $\:\tilde{}\:$ cannot be considered as a closure
operation in the usual sense, since $J\subseteq I$ does not generally
imply $\tilde J\subseteq\tilde I$. An example from \cite{RS} shows this: let
$J=\langle y^4,xy^3,x^3y,x^4\rangle\subset  I=\langle
y^3,x^3\rangle\subset k[x,y]$, then $I$ is Ratliff-Rush but
$x^2y^2\in\tilde J\backslash\tilde I$.

Several
results about Ratliff-Rush ideals are given in \cite{HJLS}, \cite{HLS}
and \cite{RR}. In addition to general results, one can find many
examples and counterexamples with respect to different properties
\cite{RS}. In \cite{E} the author presents an algorithm for computing
Ratliff-Rush associated ideals by computing the Poincar$\acute{\rm e}$
series and choosing a tame superficial sequence of $I$.

One of the reasons to study Ratliff-Rush ideals is the following. Let
$I$ be a regular $\mathfrak m$-primary ideal in a local ring
$(R,\mathfrak m,k)$. We know that the Hilbert function $H_I(l)={\rm
  dim}_k (R/I^l)$ is a polynomial $P_I(l)$ called the Hilbert polynomial
of $I$ for all large $l$. Then $\tilde I$ can be defined as the unique
largest ideal containing $I$ and having the same Hilbert polynomial as
$I$.

Ratliff-Rush ideals associated to monomial ideals are monomial by
definition, which makes the computations easier. There is always a
positive integer $L$ such that $\tilde I=I^{L+1}:I^L$, but it is not
clear how big that $L$ is (see Example~1.8 in \cite{RS}). If $I$ is a monomial
ideal and $m$ is some monomial, then for all $l\geq 0$ we have \begin{equation} \label{mIekv}
  (mI)^{l+1}:(mI)^l=(m^{l+1}I^{l+1}):(m^lI^l)=m(I^{l+1}:I^l).
\end{equation} Principal ideals are trivially Ratliff-Rush. Any
non-principal monomial ideal $J$ in the rings $k[x,y]$ and
$k[[x,y]]$ can be written as $J=mI$, where $m$ is a monomial and $I$
is an $\langle x,y\rangle$-primary ideal; hence it suffices
to consider $\langle x,y\rangle$-primary monomial
ideals. Moreover, (\ref{mIekv}) shows that the Ratliff-Rush reduction
numbers of $I$ and $mI$ are the same.

In this paper we show how to compute the Ratliff-Rush ideal
associated to a monomial ideal in a certain class in the rings $k[x,y]$ and $k[[x,y]]$ and find an upper bound
for the Ratliff-Rash reduction number for such an ideal. Section~\ref{Srons} is
devoted to some results about numerical semigroups that are crucial
for our work in Section~\ref{RRiatcmi}. In Section~\ref{E} we duscuss several useful examples.

\section{Some results on numerical semigroups} \label{Srons}

A numerical semigroup $S$ is a set of linear
combinations $\lambda_1a_1+\cdots +\lambda_ra_r$, where $a_i\in\mathbb
Z_{\geq 0}$ are the
generators and $\lambda_i\in\mathbb Z_{\geq 0}$ are the
coefficients. There is a partial ordering $\leq_{S}$ where for any pair
$s,s'$ in $S$, if there is $s''\in S$ such that $s'=s+s''$ then $s\leq
s'$. The set of minimal elements in $S\backslash\{0\}$ in this
ordering is called a $minimal\; set\; of\; generators$ for $S$. If a
semigroup is generated by a set $\{a_i\}_{i=1}^r$, then we denote it
by $\langle a_1,\ldots ,a_r\rangle$.

\begin{definition} \label{DefFrob}
Let $S=\langle a_i\rangle$ be a numerical semigroup and $\gcd
(a_i)=h$. The greatest multiple of $h$ that does not belong to $S$ is
called the $Frobenius\; number$ of $S$ and is denoted by $g(S)$. If $\gcd
(a_i)=1$, then the Frobenius number is the greatest integer that does
not belong to $S$. A list of references to the papers written about this
subject can be found in \cite{FGH}, pp. 1-2.
\end{definition}

We notice that for any $h\in\mathbb Z_+$ the numerical semigroups $\langle
a_i\rangle$ and $\langle ha_i\rangle$ are isomorphic.

\begin{definition} \label{DefLambda}
Let $S=\langle a_1,\ldots, a_r\rangle$, where $a_1<\cdots <a_r$, be a numerical semigroup. For $s\in S$ the
coefficients in a linear combination $s=\sum\lambda_ia_i$ are not
necessarily unique. We define the function $\lambda :
S\rightarrow\mathbb Z_{\geq 0}$ by $\lambda (s)=\min\,\{\,\sum\lambda_i\:\vert\:
s=\sum\lambda_ia_i\,\}$. Then we define the following positive number:
\begin{equation} \label{EqLambda} \Lambda=\Lambda(S) =\max\,\{\,\lambda(s)\:\vert\: s\leq g(S)+a_r\,\}.\end{equation}
\end{definition}

\begin{corollary}
Let $S=\langle a_1,\ldots ,a_r\rangle$ with $a_1<\cdots <a_r$. Then for $s\in S$ we have
$\lim_{s\to\infty}\frac{s}{\lambda(s)}=a_r$.
\end{corollary}

\begin{proof}
For each $s>g(S)$ there is $n\in\mathbb Z_{\geq
0}$ such that $g(S)+a_rn+1\leq s\leq g(S)+a_r(n+1).$ Then, obviously,
$\lambda(s)\geq n$ and $\lambda(s)\leq\Lambda +n$
by Definition \ref{DefLambda}. Hence, $\frac{g(S)+a_rn+1}{\Lambda
  +n}\leq\frac{s}{\lambda(s)}\leq\frac{g(S)+a_r(n+1)}{n}$. The limits of both the right hand side and the left hand
side are $a_r$ as $s\rightarrow\infty$.
\end{proof}

\begin{proposition} \label{Lemalfabeta}
Let $S=\langle a_1,\ldots ,a_r\rangle$ be a numerical semigroup
generated by nonnegative integers $a_1<\cdots <a_r$. Let $\alpha <1$ and $\beta$
be real nonnegative numbers. Then there is a number L such that for every integer $l\geq
L$ the following is true:

if $s\in S$ and $s\leq a_r\cdot \alpha l+\beta$, then $\lambda(s)\leq l$.
\end{proposition}

\begin{proof}
For each $s>g(S)$ there is $n\in\mathbb Z_{\geq 0}$ such
that $g(S)+a_rn+1\leq s\leq g(S)+a_r(n+1)$. Thus,
$\lambda(s)\leq\Lambda+n\leq\Lambda+\frac{s-g(S)-1}{a_r}$. Hence, if
$s\leq a_r\alpha l+\beta$ we get $\lambda(s)\leq\Lambda+\alpha l+\frac{\beta -g(S)-1}{a_r}$. We want to find
an $L$ such that $\lambda(s)\leq l$ for all $l\geq L$. This occures if
\begin{equation}\label{Ekvalfabeta} l\geq\frac{a_r\Lambda +\beta
    -g(S)-1}{a_r(1-\alpha)},\end{equation} which is an upper bound for
the number $L$.
\end{proof}

\begin{remark} \label{Remalfa<1}
It is easy to see the necessity of the condition $\alpha <1$. Consider the
numerical semigoup $S=\langle 2,5\rangle$. Let $l\in\mathbb Z_{\geq 0}$ and
$\beta =4$. Then for any $l$ there is no $\lambda_1\in\mathbb
Z_{\geq 0}$ such that $s=5l+4=\lambda_1\cdot 5+(l-\lambda_1)\cdot 2$.
\end{remark}

\begin{corollary} \label{Corsemigrupps+t=d}
Let $S=\langle a_i\rangle_{i=0}^r$ and $T=\langle b_i\rangle_{i=0}^r$,
where $a_0=b_r=0$ and $a_i+b_i=d$ for all $i$, be numerical semigroups. Then there is a number $L$ such that for every integer $l\geq
L$ and some fixed $\beta$ the following is true:

if $s\in S$ and $s\leq d\cdot\alpha l+\beta\leq dl$, then there are
$\lambda_0,\ldots ,\lambda_r$ such that $s=\sum_{i=0}^r\lambda_ia_i$
and $\sum_{i=0}^r\lambda_i=l$; moreover,
$dl-s=\sum_{i=0}^r\lambda_ib_i\in T$.
\end{corollary}

\begin{proof}
By Proposition~\ref{Lemalfabeta} there is some $L$ such that for all $l\geq
L$ if $S\ni s\leq d\alpha l+\beta$ then $s=\sum_{i=1}^r\lambda_ia_i$, where
$\sum_{i=1}^r\lambda_i\leq l$. Letting $\lambda_0
=l-\sum_{i=1}^r\lambda_i$ we can write $s=\sum_{i=0}^r\lambda_ia_i$ where $\sum_{i=0}^r\lambda_i=l$. Clearly, $dl-s=d\sum_{i=0}^r\lambda_i-\sum_{i=0}^r\lambda_ia_i=\sum_{i=0}^r\lambda_ib_i\in T$.
\end{proof}

If $S$ and $T$ are as in Corollary~\ref{Corsemigrupps+t=d}, then for
any $s\in S$ and $t\in T$ such that $s+t=dl$ we have either
$s\leq\frac{dl}{2}+\beta$ or $t\leq\frac{dl}{2}+\beta$ for some
$\beta\leq\frac{dl}{2}$. This estimation will be used frequently in the
next two sections when we apply our results on calculating powers of
and Ratliff-Rush ideals associated to some monomial ideals.

\begin{example} \label{Exs<dl/2}
In Proposition~\ref{Lemalfabeta} let $\alpha=\frac{1}{2}$ and
$\beta=0$. Thus, for every $l\geq 2\Lambda -\frac{2g(S)+2}{a_r}$, if
$s\in S$ and $s\leq \frac{a_rl}{2}$, then $\lambda(s)\leq l$.
\end{example}

\begin{example} \label{Exs<(dl-1)/2}
For any $l\in\mathbb Z_{\geq 0}$ every $s\in S$ belongs to the
interval $\frac{a_r(l+j)}{2}\leq s\leq\frac{a_r(l+j+1)-1}{2}$ for some
$j\geq -l$. That is, the assumptions in Proposition~\ref{Lemalfabeta}
are fulfilled for $\alpha=\frac{1}{2}$ and
$\beta=\frac{a_r(j+1)-1}{2}$. Hence, for all $l\geq 2\Lambda
+(j+1)-\frac{3+2g(S)}{a_r}$ if
$s\leq\frac{a_rl}{2}+\frac{a_r(j+1)-1}{2}$ then $\lambda(s)\leq l$.
\end{example}

\section{Ratliff-Rush ideals associated to certain monomial ideals}\label{RRiatcmi}

Now we will apply the results from the previous section in order to
compute Ratliff-Rush ideals for some monomial cases. We start with the case where all the minimal generators
of the ideal have the same degree.

\subsection{Ideals generated by monomials of the same degree} \label{Seca+b=d}

Let $I=\langle x^{a_i}y^{b_i}\rangle_{i=0}^r$ be an $\mathfrak m$-primary ideal
generated by the monomials of the same degree $d$ ordered in such a way
that $a_i<a_{i+1}$ and $b_i>b_{i+1}$; in other words, $a_0=b_r=0$ and
$b_i=d-a_i$ for all $i$. To this ideal we associate the numerical
semigroups $S=\langle a_i\rangle_{i=0}^r$ and $T=\langle b_i\rangle_{i=0}^r$.

The ideal $I^l$ is generated by monomials of degree $dl$, namely by
\begin{equation} \label{Ilgeneratorer} \{\prod_{\sum l_i=l} (x^{a_i}y^{b_i})^{l_i}=x^{\sum l_ia_i}y^{\sum
    l_ib_i}\}. \end{equation} Here $\sum l_ia_i\in S,\; \sum
l_ib_i\in T$ and $\sum l_ia_i +\sum l_ib_i=\sum l_i(a_i+b_i)=dl.$

\begin{theorem}\label{Kvadratlinjeidealsats}
Let an ideal $I=\langle x^{a_i}y^{b_i}\rangle_{i=0}^r\subset R$ and the
corresponding numerical semigroups $S=\langle a_i\rangle_{i=0}^r$ and $T=\langle
b_i\rangle_{i=0}^r$. Then there is an integer L such that for
any $l\geq L$ the following is true:
\begin{equation} \label{Ekvkvadratsats}I^l=\langle x^sy^t\:\vert\: s\in S\; {\rm and}\; t\in T\; {\rm such\; that}\;
s+t=dl\rangle .\end{equation} Moreover, for $l$ sufficiently large:
\begin{enumerate}
\item\label{su} if $s\in S,\; s\leq u$ and $s+u\geq dl$, where $u\in\mathbb
Z_{\geq 0}$, then $x^sy^u\in I^l$;

\item\label{vt} if $t\in T,\; t\leq v$ and $t+v\geq dl$, where $v\in\mathbb
Z_{\geq 0}$, then $x^vy^t\in I^l$.

\end{enumerate}
\end{theorem}

\begin{proof}
The inclusion $I^l\subseteq\langle x^sy^t\:\vert\: s\in S,\; t\in T\; {\rm and}\;
s+t=dl\rangle$ is true for all $l$, which is clear from the text preceeding the theorem.

\smallskip
The other inclusion needs to be proved since $s+t=dl$ does not
generally imply that $s=\sum_{i=0}^r\lambda_ia_i$ with
$\sum_{i=0}^r\lambda_i=l$ or $t=\sum_{i=0}^r\mu_ib_i$ with $\sum_{i=0}^r\mu_i=l$. However, this is
asserted by Corollary~\ref{Corsemigrupps+t=d} as we will see
below. Thus, we will prove the second part of the theorem, because this other
inclusion is a special case of it, since if $s+t=dl$ then either
$s\leq t$ or
$t\leq s$.

\medskip
(\ref{su}) If $s+u\geq dl$ then there is some $j$ such that
$d(l+j)\leq s+u\leq d(l+j+1)-1$. We will show that for any such $s,u$
and $j$ we have $x^sy^u\in I^{l+j}\subset I^l$. Clearly, it is
sufficient to consider the case $j=0$, that is, suppose $dl\leq
s+u=dl+b\leq d(l+1)-1$. If $s\leq u$ then
$s\leq\frac{dl}{2}+\frac{d-1}{2}$. By Corollary~\ref{Corsemigrupps+t=d}, for
sufficiently large $l$ we can write $s=\sum_{i=0}^r\lambda_ia_i$ and
$u=b+\sum_{i=0}^r\lambda_ib_i$, where
$\sum\lambda_i=l$. Hence, $x^sy^u=y^b\prod_i
(x^{a_i}y^{b_i})^{\lambda_i}\in I^l$.

\smallskip
Part (\ref{vt}) is proved similarly.
\end{proof}

\begin{remark} \label{Kvadratpotenslimsup}
By Example~\ref{Exs<(dl-1)/2} an upper bound for the least integer $L$
in Theorem \ref{Kvadratlinjeidealsats} is $\lceil\max\big(
2\Lambda(S)+1-\frac{3+2g(S)}{d},2\Lambda(T)+1-\frac{3+2g(T)}{d}\big)\rceil$,
where $\lceil c\rceil$ denotes the least integer which is greater or
equal to $c$.
\end{remark}

\begin{definition} \label{DefISIT}
Let the assumptions be as in Theorem~\ref{Kvadratlinjeidealsats}. We introduce
the following ideals: $I_S=\langle x^sy^{d-s}\:\vert\: s\in S\; {\rm and}\; s\leq
d\rangle$ and $I_T=\langle x^{d-t}y^t\:\vert\: t\in T\; {\rm and}\; t\leq d\rangle$.
\end{definition}

\begin{proposition} \label{PropISIT}
Let the assumptions be as in Theorem~\ref{Kvadratlinjeidealsats}. Then
for every l sufficiently large \begin{equation} \label{Ilsomtresumma}
  I^l=\big(y^{dl-d}\big)I_S+\big(x^{dl-d}\big)I_T+\big(x^dy^d\big)I_{M,l}\end{equation} for some ideal $I_{M,l}$.
\end{proposition}

\begin{proof}
Let $l\geq\max(\{\lambda(s)\:\vert\: S\ni s\leq d
\},\{\lambda(t)\:\vert\: T\ni t\leq d \})$.

Then $y^{dl-d}(x^sy^{d-s})\in I^l$ if and only if
$s=\sum\lambda_i a_i\leq d$ where $\sum\lambda_ia_i=l$. Equivalently,
$x^{dl-d}(x^{d-t}y^t)\in I^l$ if and only if $t=\sum\lambda_ib_i\leq
d$. Finally, the generators for $I^l$ such that both the power of $x$
and $y$ is equal to or greater than $d$ can be written as the third term in (\ref{Ilsomtresumma}) where $I_{M,l}=I:(x^dy^d)$.
\end{proof}

\begin{example} \label{ExIlsomsumma}
Let $I=\langle y^7,x^2y^5,x^5y^2,x^7\rangle$. Then $I_S=\langle
y^7,x^2y^5,x^4y^3, x^5y^2,\allowbreak x^6y,x^7\rangle$ and $I_T=\langle
y^7,xy^6,x^2y^5,x^3y^4,x^5y^2,x^7\rangle$. For $l\geq 3$ we can write
$I^l=y^{7l-7}I_S+x^{7l-7}I_T+x^7y^7I_{M,l}$ for some $I_{M,l}$. For $l\geq 4$
the ideal $I_{M,l}=\mathfrak m^{7l-14}$.
\end{example}

\begin{remark}\label{IMlmaxpotens}
Generally, if $g(S)$ and $g(T)$ are less or equal to $d-1$, then
$I_{M,l}=\mathfrak m^{d(l-2)}$ for all sufficiently large $l$.
\end{remark}

\begin{proposition} \label{KvadratRR}
Let $I=\langle x^{a_i}y^{b_i}\rangle_{i=0}^r\subset R,\; S=\langle a_i\rangle_{i=0}^r$ and
$T=\langle b_i\rangle_{i=0}^r$.

Then the Ratliff-Rush ideal associated to I is $$\tilde I=I_S\cap I_T .$$
\end{proposition}

\begin{proof}
We will show that $I^{l+1}:I^l=I_S\cap I_T$ for all sufficiently large
$l$. Since $I$ is monomial, a polynomial $p$ belongs to $I^{l+1}$ if and only if
every power product in $p$ belongs to $I^{l+1}$. Hence, it suffices to
consider monomial ring elements.

\medskip
Let $m\in I_S\cap I_T$. Then
$m=m'x^{s'}y^{d-s'}=m''x^{d-t''}y^{t''}$. We know that for all
sufficiently large $l$ the generators
for $I^l$ are on the form $x^sy^t$ where $s\in S,\; t\in T$ and
$s+t=dl$, that is either $s\leq\frac{dl}{2}$ or $t\leq\frac{dl}{2}$. Assume $s\leq\frac{dl}{2}$. Then, using the first equality
for $m$, we get $m\cdot x^sy^t=m'x^{s+s'}y^{dl+d-(s+s')}$. Since
$s+s'\leq\frac{dl}{2}+d=\frac{d(l+1)}{2}+\frac{d}{2}$, then by Corollary~\ref{Corsemigrupps+t=d} there is some
integer $L_S$ use we can write such that for all $l\geq L_S$ we can write
$s+s'=\sum\lambda_ia_i$ and $d(l+1)-(s+s')=\sum\lambda_ib_i$ where $\sum\lambda_i
=l+1$. Hence, $mx^sy^t=\prod_i (x^{a_i}y^{b_i})^{\lambda_i}\in
I^{l+1}$.

Using the equality $m=m''x^{d-t''}y^{t''}$ and Corollary~\ref{Corsemigrupps+t=d} we show in the same way that there is some $L_T$ such that for all $l+1\geq L_T$ if
$t\leq\frac{dl}{2}$ then $mx^sy^t=\prod_i
(x^{a_i}y^{b_i})^{\mu_i}\in I^{l+1}$.

\medskip
On the other hand, assume $m\notin I_S$. Then $my^{dl}\notin
y^{dl}I_S$ and, hence, $my^{dl}\notin I^{l+1}$ by Proposition 
\ref{PropISIT}. Analogously, if $m\notin I_T$ then $mx^{dl}\notin
I^{l+1}$, which finishes the proof.
\end{proof}

\begin{corollary} \label{S=a_i}
Let $I=\langle x^{a_i}y^{b_i}\rangle_{i=0}^r\subset R,\; S=\langle a_i\rangle_{i=0}^r$ and
$T=\langle b_i\rangle_{i=0}^r$ its corresponding numerical semigroups. If for every pair $a_i$ and
$a_j$ we have either $a_i+a_j=a_k$ for some k or $a_i+a_j\geq d$, then
$I$ is Ratliff-Rush.
\end{corollary}

\begin{proof}
Clearly, the set $\{s\in S\:\vert\: s\leq d\}=\{a_i\}_{i=0}^r$ and then
$I_S=I$. Since the inclusion $I\subseteq I_T$ is always valid, we
conclude that $\tilde I=I_S\cap I_T=I\cap I_T=I$. 
\end{proof}

\begin{proposition}\label{KvadratRRlimsup}
Let I, S and T be as in Theorem~\ref{Kvadratlinjeidealsats}. Then there is an upper bound for the reduction number
of I: \begin{equation} r(I)\leq\lceil\max\big(
  2\Lambda(S)+2-\frac{g(S)+1}{d},2\Lambda(T)+2-\frac{g(T)+1}{d}\big)\rceil -1.\end{equation}
\end{proposition}

\begin{proof} The proof of Proposition~\ref{KvadratRR} asserts that
  the upper bound is, using the notations from there, equal to
  $\lceil\max(L_S,L_T)\rceil$. The result follows from the formula (\ref{Ekvalfabeta}) in Proposition~\ref{Lemalfabeta} with $\alpha =\frac{1}{2}$ and $\beta =d$ .
\end{proof}

\begin{example} 
Let $I$ be the ideal in Example~\ref{ExIlsomsumma}. Then $\tilde I=I_S\cap I_T
=\langle y^7,x^2y^5, x^4y^4,x^5y^2,x^7\rangle$. It is
interesting to note that $I^l$ satisfies (\ref{Ekvkvadratsats})
for all $l\geq 5$ by Remark~\ref{Kvadratpotenslimsup}, but actually
for all $l\geq 4$. Further, $r(I)=1$ while the upper bound suggested
by Proposition~\ref{KvadratRRlimsup} is five.
\end{example}

\begin{example}
Let $I=\langle y^{18}, x^3y^{15}, x^{13}y^5, x^{18}\rangle$. Then
$\tilde I=I_S\cap I_T=\langle y^{18},\allowbreak
  x^3y^{15}, x^8y^{12},\allowbreak x^9y^{10}, x^{13}y^5, x^{18}\rangle$
  and $r(I)=4$. Thus, the minimal generators for $\tilde I$ do not
  need to be of the same degree.
\end{example}

\subsection{Ideals generated by $x^{a_i}y^{b_i}$ such that $\frac{a_i}{a_r}+\frac{b_i}{b_0}=1$}

Here we discuss slight generalizations of the subject in Section~\ref{Seca+b=d} to $\langle x,y\rangle$-primary monomial ideals $\langle
x^{a_i}y^{b_i}\rangle_{i=0}^r$ such that $\frac{a_i}{a_r}+\frac{b_i}{b_0}=1$ where $\gcd(b_0,a_r)=d$. We
can, of course, apply the results directly using the numerical
semigroups $S'=\frac{d}{a_r}\cdot S$ and $T'=\frac{d}{b_0}\cdot
T$. However, it might be useful to devote some space to formulate the
material differently in order to make it possible to widen
the results.

\begin{corollary} \label{Snedlinjesemigrupp}
Let $S=\langle a_i\rangle_{i=0}^r$ and $T=\langle b_i\rangle_{i=0}^r$,
where $a_0=b_r=0$ and $\frac{a_i}{a_r}+\frac{b_i}{b_0}=1$ for all $i$, be numerical semigroups. Then there is a number $L$ such that for every integer $l\geq
L$ and some fixed $\beta$ the following is true:

if $s\in S$ and $s\leq a_r\cdot\alpha l+\beta\leq dl$ then there are
$\lambda_0,\ldots ,\lambda_r$ such that $s=\sum_{i=0}^r\lambda_ia_i$
and $\sum_{i=0}^r\lambda_i=l$; moreover,
$l-\frac{s}{a_r}=\frac{1}{b_0}\sum_{i=0}^r\lambda_ib_i\in \frac{1}{b_0}\cdot T$.
\end{corollary}

\begin{proof} 
The proof differs from the one of Corollary~\ref{Corsemigrupps+t=d} by the last sentence, which here should be:

$b_0(l-\frac{s}{a_r})=b_0(\sum\lambda_i
-\frac{1}{a_r}\sum\lambda_ia_i)=b_0\big(\sum\lambda_i(1-\frac{a_i}{a_r})\big)=\sum\lambda_ib_i\in
T$.
\end{proof}

\begin{theorem} \label{Snedlinjepotenssats}
Let $I=\langle x^{a_i}y^{b_i}\rangle_{i=0}^r\subset R$ be an $\mathfrak
m$-primary ideal such that $a_0=b_r=0$ and
$\frac{a_i}{a_r}+\frac{b_i}{b_0}=1$. Let $S=\langle
a_i\rangle_{i=0}^r$ and $T=\langle b_i\rangle_{i=0}^r$ be numerical semigroups. Then there is an integer L such that for
any $l\geq L$ the following is true: \begin{equation}
  \label{Snedlinjepotensekv} I^l=\langle x^sy^t\:\vert\: s\in S\; {\rm and}\; t\in T\; {\rm
  such\; that}\; \frac{s}{a_r}+\frac{t}{b_0}=l
\rangle .\end{equation}

Moreover, for $l$ sufficiently large:
\begin{enumerate}
\item if $s\in S,\; \frac{s}{a_r}\leq\frac{u}{b_0}$ and
$\frac{s}{a_r}+\frac{u}{b_0}\geq l$ for some $u\in\mathbb
Z_{\geq 0}$, then $x^sy^u\in I^l$;

\item if $t\in T,\; \frac{t}{b_0}\leq\frac{v}{a_r}$ and
$\frac{t}{b_0}+\frac{v}{a_r}\geq l$ for some $v\in\mathbb
Z_{\geq 0}$, then $x^vy^t\in I^l$.
\end{enumerate}
\end{theorem}

\begin{proof}
The ideal $I^l$ is a subideal of the right hand side of
(\ref{Snedlinjepotensekv}) by the definition of $S$ and $T$ and the condition on
the exponents.

\medskip
To prove (\ref{su}) it suffices to show that if
$l\leq\frac{s}{a_r}+\frac{u}{b_0}=l+q\leq l+1$ for some rational $q$ then $x^sy^u\in I^l$;
compare to the proof of Theorem~\ref{Kvadratlinjeidealsats}.

If $\frac{s}{a_r}\leq\frac{u}{b_0}$ and
$\frac{s}{a_r}+\frac{u}{b_0}\leq l+1$, then $s\leq\frac{a_rl}{2}+\frac{a_r}{2}$. Thus, by Corollary~\ref{Snedlinjesemigrupp}, for
sufficiently large $l$ we can write $s=\sum\lambda_ia_i$ where
$\sum\lambda_i=l$. Further, let some $t=\sum\lambda_ib_i$, then
$\frac{s}{a_r}+\frac{t}{b_0}=\sum\lambda_i(\frac{a_i}{a_r}+\frac{b_i}{b_0})=l$.
Hence, $u\geq t$ and we
get $x^sy^u\in I^l$.

\smallskip
Part (\ref{vt}) is proved similarly.
\end{proof}

\begin{proposition} \label{SnedRR}
Let $I=\langle x^{a_i}y^{b_i}\rangle\subset R,\; S$ and
$T$ be as in Theorem~\ref{Snedlinjepotenssats}. Then the Ratliff-Rush ideal associated to I is \begin{equation}
  \begin{split} \tilde
  I=& \langle x^sy^u\:\vert\: s\in S,\; s\leq
  a_r\;{\rm and}\; u\;{\rm such\;
    that}\;\frac{s}{a_r}+\frac{u}{b_0}=1\rangle\cap\\
    & \langle x^vy^t\:\vert\: t\in T,\; t\leq
  b_0\;{\rm and}\; v\;{\rm such\;
    that}\;\frac{v}{a_r}+\frac{t}{b_0}=1\rangle .\end{split} \end{equation}
\end{proposition}

\begin{example}
Let $I=\langle y^{12},x^{6}y^{8},x^{9}y^{6},x^{15}y^{2},x^{18}
\rangle$. Then we have $S=\langle 6,9\rangle,\; T=\langle 2\rangle$ and $I_S=\langle y^{12},x^{6}y^{8},x^{9}y^{6},
x^{12}y^{4},x^{15}y^{2},x^{18}\rangle ,\; I_T=\langle
y^{12},x^{3}y^{10},\allowbreak x^{6}y^{8}, x^{9}y^{6},
x^{12}y^{4},x^{15}y^{2},x^{18}\rangle$. Thus, the Ratliff-Rush
associated ideal is $\tilde I=I_S\cap I_T =I+\langle x^{12}y^4\rangle$.
\end{example}

\section{Examples}\label{E}

In the sequel we let $I=\langle
x^{a_i}y^{b_i}\rangle_{i=0}^r$ be an $\langle x,y\rangle$-primary
ideal such that $a_i+b_i=d$ for all $i$.

\begin{example} \label{a_i<d/2}
Assume $a_1\geq\frac{d}{2}$, then $a_i+a_j\geq d$ for all pairs of
$a_i$ and $a_j$. Hence, the condition in Corollary~\ref{S=a_i} is
fulfilled and $I$ is Ratliff-Rush. This generalizes
the example of a non integrally closed Ratliff-Rush ideal $\langle
y^4,x^2y^2,x^3y,x^4\rangle$ in \cite{RS}, p.~2.

The only integrally
closed monomial ideals such that the generators have the same degree
are $\langle x,y\rangle^d$.

By the total
$N=2^{d-1}$ of $\langle x,y\rangle$-primary ideals generated by degree
$d$ monomials there are $2\cdot 2^{\lceil\frac{d}{2}\rceil}$ such that $\frac{d}{2}\leq a_1$
or $\frac{d}{2}\leq b_1$. Hence, such monomial ideals generated by the same degree there are $2\sqrt N$
Ratliff-Rush ideals if $d$ is odd and $2\sqrt{2N}$ if $d$ is even.
\end{example}

\subsection{Ideals such that all their powers are Ratliff-Rush}

It is shown in \cite{HLS}, (1.2), that all the powers of a regular ideal
in a Noetherian ring are Ratliff-Rush if and only if the depth of the
associated graded ring $gr_I(R)$ is positive.

\begin{example} \label{Ex6.3}
In \cite{HJLS}, (6.3), the authors conjecture that for any $d$ the ideal $I_d=\langle
y^d,x^{d-1}y,x^d\rangle$ and all its powers are Ratliff-Rush. The conjecture was later
proved in \cite{RS} by actual computation of the depth. An
alternative way to show this uses Corollary \ref{S=a_i}.

The numerical
semigroups associated to the ideal $I_d$ are
$S_d=\bigcup_{l=0}^{\infty} \{ ld-i\}_{i=0}^l$ and
$T_d=\mathbb Z_{\geq 0}$. Obviously, if $s\in S_d$ and $s\leq dl$, then $\lambda(s)\leq l$. Let $S_{d,l}$ be the numerical
semigroup associated to the ideal $I_d^l$. Then $\{s\in
S_{d,l}\:\vert\: s\leq dl\} =\{ {\rm exponents\; of}\; x\; {\rm in\;
  the\; minimal\; generating\; set\; for}\; I_d^l\}$. Hence, $I_d^l$ is
Ratliff-Rush for all $l$ by Corollary \ref{S=a_i}.

\medskip
This family of ideals is part of a larger family in which all the
powers of an ideal are Ratliff-Rush.

\medskip
Let $I_{d,k}=\langle y^d,x^{d-k}y^k, x^{d-k+1}y^{k-1}, \ldots
,x^{d-1}y,x^d\rangle$. For example, the family $I_{d,1}$ are the
ideals we have discussed previously. The corresponding numerical semigroups are
$S_{(d,k)}=\bigcup_{l=0}^{\infty}\{ ld-i\}_{0\leq i\leq lk}$ and $T_{(d,k)}=\mathbb
Z_{\geq 0}$. If $s\in S_{(d,k)}$ and $s\leq ld$ then $\lambda(s)\leq
l$. Then the exponents of $x$ among the generators for $I_{d,k}^l$
fulfil the assumption in Corollary~\ref{S=a_i}, which finishes the proof.
\end{example}

\begin{example}
In \cite{HJLS}, (E3), the authors examine the ideal $I=\langle
y^8, x^3y^5, \allowbreak x^5y^3, x^8\rangle$ using MACAULAY. Among other things they
show that $I$ is Ratliff-Rush but $I^3$ is not. We will look at all
the powers $I^l$.

Using Proposition \ref{KvadratRR} we see that $\tilde I=\langle
y^8,x^3y^5,x^5y^3,x^6y^2,x^8\rangle\cap\langle
y^8,x^2y^6,\allowbreak x^3y^5,x^5y^3,x^8\rangle =I$.

Further, $g(S)=g(\langle 0,3,5,8\rangle)=7$ and
$\lambda(s)\leq 4$ if $s\leq 24$. Thus, for every $s\leq 8k$ we have
$\lambda(s)\leq k+1$, that is, for all $l\geq 4$ if
$s\leq 8\cdot\frac{l}{2}$ then $\lambda(s)\leq\lceil\frac{l}{2}\rceil
+1\leq l$. Exactly the same is valid for the numerical semigroup $T$. Hence, \begin{equation}\label{E3^l} I^l=\langle
  x^sy^{8l-s}\:\vert\: s\in S\; {\rm and}\; s\leq 4l\rangle +\langle
  x^{8l-t}y^t\:\vert\: t\in T\; {\rm and}\; t\leq 4l\rangle\end{equation} for all $l\geq 4$. (Compare to Remark
\ref{IMlmaxpotens}.) Moreover,
$I^2$ is on that form too, but not $I^3$ since $\lambda(12)=4$.

Now we will show that if $I^l$ is on the form (\ref{E3^l}), then $I^l$
is Ratliff-Rush. Let $S_l$ and $T_l$ be the numerical semigroups
defined by $I^l$, then $I_{S_l}=\langle x^sy^{8l-s}\:\vert\: s\leq 4l\rangle
+x^{4l}\mathfrak m^{4l}$ and $I_{T_l}=\langle x^{8l-t}y^t\:\vert\:
t\leq 4l\rangle +y^{4l}\mathfrak m^{4l}$. It is easy to see that
$I_{S_l}\cap I_{T_l}=\langle x^sy^{8l-s}\:\vert\: s\leq 4l\rangle +\langle x^{8l-t}y^t\:\vert\:
t\leq 4l\rangle =I^l$.

Finally, we get $\widetilde{I^3}=I^3+\langle x^{12}y^{12}\rangle$ using Proposition \ref{KvadratRR}.
\end{example}

\begin{example}
Let $I_k=\langle y^{6k+1}\rangle +\left\langle
x^{2(k+i)+1}y^{4k-2i}\right\rangle_{i=0}^{k-1}+\left\langle
x^{4k+i+1}y^{2k-i}\right\rangle_{i=0}^{2k}$.

\smallskip For example, if $k=2$ then
$I_2=\langle
y^{13},x^5y^8,x^7y^6,x^9y^4,x^{10}y^3,x^{11}y^2,x^{12}y,x^{13}\rangle$.

\medskip
We will prove that all positive powers of $I_k$ are Ratliff-Rush by showing that
the numerical semigroup determined by $I$ is $S=\{a_i\}\cup\{n\in\mathbb
Z\:\vert\: n\geq 6k+1\}$ and if $s\in S$ is such that $s\leq l(6k+1)$
then $\lambda(s)\leq l$. Hence, the generators for $I_k^l$ will fulfil
the condition in Corollary \ref{S=a_i}.

We use induction on $l$.

If $l=1$ we are done, since $\{ s\in
S\:\vert\: s\leq 6k+1\} =\{a_i\}$.

Let $l=2$. We will show that all the elements in $\{ 6k+2,\ldots
,12k+2\}$ are linear combinations of at most two generators $a_i$ and
$a_j$. For all $0\leq i\leq 2k$ we have $6k+2\leq (2k+1)+(4k+i+1)\leq
8k+2$. Further, any integer $n\in [8k+2,\ldots ,12k+2]$ is a linear
combination of two elements in $\{ 4k+i+1\}_{i=0}^{2k}$.

Assume our claim is true for all
$l\leq p$. Let $l=p+1$. We need to show that if $p(6k+1)+1\leq n\leq
l(6k+1)+6k+1$ then $n=\sum \lambda_ia_i$ with $\sum\lambda_i\leq
p+1$. By the induction hypothesis $\{
p(6k+1)-4k+i\}_{i=0}^{4k}\subset S$ and the values of the
$\lambda$-function of these elements are always less or equal to
$p$. Thus, $\big(p(6k+1)-4k+i\big)+4k+1=\sum\lambda_ia_i$ with
$\sum\lambda_i\leq p+1$ for all $0\leq i\leq 4k$. Clearly, the same is
valid for each sum $p(6k+1)+(4k+i+1)$ for all $0\leq i\leq 2k$, and we
are done.

\medskip
This last example can be varied in many different ways. Moreover,
the induction proof that we used can be applied on other families of
ideals. For example, $I_{n,k}=\left\langle
x^{in}y^{n(k+1-i)-1}\right\rangle_{i=0}^k+\left\langle
x^{kn+j}y^{n-j-1}\right\rangle_{j=0}^{n-1}$. If $n=3$ we get the
family $I_{3,k}=\langle
y^{3k+2},x^3y^{3k-1},x^6y^{3k-4},\ldots ,x^{3k}y^2,x^{3k+1}y,x^{3k+2}\rangle$.
\end{example}

\section*{Acknowledgement}
I want to thank my former advisors Ralf Fr${\rm \ddot{o}}$berg and
Christian Gottlieb for introducing me to the subject and, most important, for many
helpful discussions and useful comments during the work on this paper.

\bibliographystyle{amsalpha}

\begin{thebibliography}{A}

\bibitem[C]{C} V. Crispin Qui$\tilde{\rm n}$onez, \textit{Integrally
    Closed Monomial Ideals and Powers of Ideals}, Research Reports in
  Mathematics {\bf 7}, Department of Mathematics, Stockholms universitet, 2002.

\bibitem[E]{E} J. Elias, \textit{On the Computation of the Ratliff-Rush
  Closure}, J. Symbolic Comput. {\bf 37} (2004), no. 6, 717-725.

\bibitem[FGH]{FGH} R. Fr${\rm \ddot{o}}$berg, C. Gottlieb and R. H${\rm
    \ddot{a}}$ggkvist, \textit{On Numerical
  Semigroups}, Semigroup Forum {\bf 35} (1987), 63-83.

\bibitem[HJLS]{HJLS} W. Heinzer, B. Johnston, D. Lantz and K. Shah, \textit{Coefficient Ideals in and Blowups of a Commutative Noetherian Domain}, J. of Algebra {\bf 162} (1993), 355-391.

\bibitem[HLS]{HLS} W. Heinzer, D. Lantz and K. Shah, \textit{The Ratliff-Rush
  Ideals in a Noetherian Ring}, Comm. in Algebra {\bf 20} (1992), 591-622.

\bibitem[RR]{RR} L. J. Ratliff, Jr and D. E. Rush, \textit{Two Notes on
    Reductions of Ideals}, Indiana Univ. Math. J. {\bf 27} (1978),
  no. 6, 929-934.

\bibitem [RS]{RS} M. E. Rossi and I. Swanson, \textit{Notes on the Behavior of
  the Ratliff-Rush Filtration}, Commutative Algebra
    (Grenoble/Lyon, 2001), 313-328, Contemp. Math., 331, Amer. Math. Soc., Providence, RI, 2003.

\end{thebibliography}

\end{document}